# Addendum to
# Classification of irreducible holonomies
# of torsion-free affine connections

By Sergei Merkulov and Lorenz Schwachhöfer*

The real form $\mathrm{Spin}(6,\mathbb{H}) \subset \mathrm{End}(\mathbb{R}^{32})$ of $\mathrm{Spin}(12,\mathbb{C}) \subset \mathrm{End}(\mathbb{C}^{32})$ is absolutely irreducible and thus satisfies the algebraic identities (40) and (41). Therefore, it also occurs as an exotic holonomy and the associated supermanifold $\mathcal{M}_{\mathfrak{g}}$ admits a SUSY-invariant polynomial. This real form has been erraneously omitted in our paper.

Also, the two real four-dimensional exotic holonomies, whose occurrences were unknown at the time of writing, have been shown to exist very recently by R. Bryant [B].

With these corrections, Table 3 and the table in Theorem C should read as follows.

---

**Table 3:** List of exotic holonomies

| group $G$ | representation $V$ | restrictions/remarks |
|---|---|---|
| $T_{\mathbb{R}} \cdot \mathrm{Spin}(5,5)$ | $\mathbb{R}^{16}$ | |
| $T_{\mathbb{R}} \cdot \mathrm{Spin}(1,9)$ | $\mathbb{R}^{16}$ | |
| $T_{\mathbb{C}} \cdot \mathrm{Spin}(10,\mathbb{C})$ | $\mathbb{C}^{16} \simeq \mathbb{R}^{32}$ | |
| $T_{\mathbb{R}} \cdot E_6^1$ | $\mathbb{R}^{27}$ | |
| $T_{\mathbb{R}} \cdot E_6^4$ | $\mathbb{R}^{27}$ | |
| $T_{\mathbb{C}} \cdot E_6^{\mathbb{C}}$ | $\mathbb{C}^{27} \simeq \mathbb{R}^{54}$ | |
| $T_{\mathbb{R}} \cdot \mathrm{SL}(2,\mathbb{R})$ | $\odot^3 \mathbb{R}^2 \simeq \mathbb{R}^4$ | |
| $\mathrm{SL}(2,\mathbb{C})$ | $\odot^3 \mathbb{C}^2 \simeq \mathbb{R}^8$ | |
| $\mathbb{C}^* \cdot \mathrm{SL}(2,\mathbb{C})$ | $\odot^3 \mathbb{C}^2 \simeq \mathbb{R}^8$ | |
| $\mathbb{R}^* \cdot \mathrm{Sp}(2,\mathbb{R})$ | $\mathbb{R}^4$ | |
| $\mathbb{C}^* \cdot \mathrm{Sp}(2,\mathbb{C})$ | $\mathbb{C}^4 \simeq \mathbb{R}^8$ | |
| $\mathbb{R}^* \cdot \mathrm{SO}(2) \cdot \mathrm{SL}(2,\mathbb{R})$ | $\mathbb{R}^2 \otimes \mathbb{R}^2 \simeq \mathbb{R}^4$ | |
| $\mathbb{C}^* \cdot \mathrm{SU}(2)$ | $\mathbb{C}^2 \simeq \mathbb{R}^4$ | |
| $H_\lambda \cdot \mathrm{SU}(2)$ | $\mathbb{C}^2 \simeq \mathbb{R}^4$ | |
| $H_\lambda \cdot \mathrm{SU}(1,1)$ | $\mathbb{C}^2 \simeq \mathbb{R}^4$ | |
| $\mathrm{SL}(2,\mathbb{R}) \cdot \mathrm{SO}(p,q)$ | $\mathbb{R}^2 \otimes \mathbb{R}^{p+q} \simeq \mathbb{R}^{2(p+q)}$ | $p+q \geqslant 3$ |
| $\mathrm{Sp}(1) \cdot \mathrm{SO}(n,\mathbb{H})$ | $\mathbb{H}^n \simeq \mathbb{R}^{4n}$ | $n \geqslant 2$ |
| $\mathrm{SL}(2,\mathbb{C}) \cdot \mathrm{SO}(n,\mathbb{C})$ | $\mathbb{C}^2 \otimes \mathbb{C}^n \simeq \mathbb{R}^{4n}$ | $n \geqslant 3$ |
| $E_7^5$ | $\mathbb{R}^{56}$ | |
| $E_7^7$ | $\mathbb{R}^{56}$ | |
| $E_7^{\mathbb{C}}$ | $\mathbb{R}^{112} \simeq \mathbb{C}^{56}$ | |
| $\mathrm{Sp}(3,\mathbb{R})$ | $\mathbb{R}^{14} \subset \Lambda^3 \mathbb{R}^6$ | |
| $\mathrm{Sp}(3,\mathbb{C})$ | $\mathbb{R}^{28} \simeq \mathbb{C}^{14} \subset \Lambda^3 \mathbb{C}^6$ | |
| $\mathrm{SL}(6,\mathbb{R})$ | $\mathbb{R}^{20} \simeq \Lambda^3 \mathbb{R}^6$ | |
| $\mathrm{SU}(1,5)$ | $\mathbb{R}^{20}$ | |
| $\mathrm{SU}(3,3)$ | $\mathbb{R}^{20}$ | |
| $\mathrm{SL}(6,\mathbb{C})$ | $\mathbb{R}^{40} \simeq \Lambda^3 \mathbb{C}^6$ | |
| $\mathrm{Spin}(2,10)$ | $\mathbb{R}^{32}$ | |
| $\mathrm{Spin}(6,6)$ | $\mathbb{R}^{32}$ | |
| $\mathrm{Spin}(6,\mathbb{H})$ | $\mathbb{R}^{32}$ | |
| $\mathrm{Spin}(12,\mathbb{C})$ | $\mathbb{C}^{32} \simeq \mathbb{R}^{64}$ | |

Notation: $T_{\mathbb{F}}$ denotes any connected Lie subgroup of $\mathbb{F}^*$, $H_\lambda = \left\{ e^{(2\pi i + \lambda)t} \mid t \in \mathbb{R} \right\} \subseteq \mathbb{C}^*$, $\lambda > 0$.

## Table from Theorem C

| Group $G$ | Representation space | Group $G$ | Representation space |
|---|---|---|---|
| $\mathrm{Sp}(n,\mathbb{R})$ | $\mathbb{R}^{2n}$ | $\mathrm{E}_7^5$ | $\mathbb{R}^{56}$ |
| $\mathrm{Sp}(n,\mathbb{C})$ | $\mathbb{C}^{2n}$ | $\mathrm{E}_7^7$ | $\mathbb{R}^{56}$ |
| $\mathrm{SL}(2,\mathbb{R})$ | $\mathbb{R}^4 \simeq \odot^3 \mathbb{R}^2$ | $\mathrm{E}_7^{\mathbb{C}}$ | $\mathbb{C}^{56}$ |
| $\mathrm{SL}(2,\mathbb{C})$ | $\mathbb{C}^4 \simeq \odot^3 \mathbb{C}^2$ | $\mathrm{Spin}(2,10)$ | $\mathbb{R}^{32}$ |
| $\mathrm{SL}(2,\mathbb{R}) \cdot \mathrm{SO}(p,q)$ | $\mathbb{R}^{2(p+q)}$, $p+q \geqslant 3$ | $\mathrm{Spin}(6,6)$ | $\mathbb{R}^{32}$ |
| $\mathrm{SL}(2,\mathbb{C}) \cdot \mathrm{SO}(n,\mathbb{C})$ | $\mathbb{C}^{2n}$, $n \geqslant 3$ | $\mathrm{Spin}(6,\mathbb{H})$ | $\mathbb{R}^{32}$ |
| $\mathrm{Sp}(1)\mathrm{SO}(n,\mathbb{H})$ | $\mathbb{H}^n \simeq \mathbb{R}^{4n}$, $n \geqslant 2$ | $\mathrm{Spin}(12,\mathbb{C})$ | $\mathbb{C}^{32}$ |
| $\mathrm{SL}(6,\mathbb{R})$ | $\mathbb{R}^{20} \simeq \Lambda^3 \mathbb{R}^6$ | $\mathrm{Sp}(3,\mathbb{R})$ | $\mathbb{R}^{14} \subset \Lambda^3 \mathbb{R}^6$ |
| $\mathrm{SU}(1,5)$ | $\mathbb{R}^{20}$ | $\mathrm{Sp}(3,\mathbb{C})$ | $\mathbb{C}^{14} \subset \Lambda^3 \mathbb{C}^6$ |
| $\mathrm{SU}(3,3)$ | $\mathbb{R}^{20}$ | | |
| $\mathrm{SL}(6,\mathbb{C})$ | $\mathbb{C}^{20} \simeq \Lambda^3 \mathbb{C}^6$ | | |

GLASGOW UNIVERSITY, GLASGOW, UK
*E-mail address*: sm@maths.gla.ac.uk

MATHEMATISCHES INSTITUT, UNIVERSITÄT LEIPZIG, LEIPZIG, GERMANY
*E-mail address*: schwachh@mathematik.uni-leipzig.de